\documentclass[letterpaper, 10 pt, conference]{ieeeconf}  

\IEEEoverridecommandlockouts                              
\usepackage{amsfonts}
\usepackage{graphicx}
\usepackage[dvipsnames]{xcolor}
\usepackage{amsmath}
\usepackage{amssymb}
\usepackage{dsfont}

\usepackage{ifthen}
\usepackage{color}
\usepackage{cite}

\usepackage{mathtools}
\usepackage{booktabs}
\usepackage{url}
\usepackage{hyperref}
\usepackage{caption}
\usepackage{subfig}
\usepackage{float}
\usepackage{multirow}
\usepackage{comment}
\usepackage{tikz}
\usepackage{pgfplots}
\pgfplotsset{compat=1.12}

\def\scr#1{{\cal #1}}
\newcommand{\R}{{\rm I\!R}}
\def\eq#1{\begin{equation}#1\end{equation}}

\newcommand{\bbb}{\mathbb}
\newtheorem{theorem}{Theorem}
\newtheorem{lemma}{Lemma}

\newtheorem{definition}{Definition}

\newtheorem{corollary}{Corollary}

\def\qed{ \rule{.08in}{.08in}}

\newcommand{\1}{\mathbf{1}}

\title{\LARGE \bf Fast Consensus Topology Design via Minimizing Laplacian Energy
}

\author{Susie Lu \hspace{.3in} Ji Liu
\thanks{
J. Liu is with the Department of Electrical and Computer Engineering at Stony Brook University
(\texttt{ji.liu@stonybrook.edu}).
S.~Lu is with Stanford Online High School (\texttt{SusieLu@ohs.stanford.edu}).
}
}

\begin{document}

\maketitle
\thispagestyle{empty}
\pagestyle{empty}



\begin{abstract}
This paper characterizes the graphical properties of an optimal topology with minimal Laplacian energy under the constraint of fixed numbers of vertices and edges, and devises an algorithm to construct such connected optimal graphs. These constructed graphs possess maximum vertex and edge connectivity, and more importantly, exhibit large algebraic connectivity of an optimal order provided they are not sparse. These properties guarantee fast and resilient consensus processes over these graphs.
\end{abstract}

\section{Introduction}

Over the past two decades, consensus has achieved great success and attracted significant attention \cite{vicsekmodel,reza1,luc,ReBe05}, being applied to a wide range of distributed control and computation problems \cite{fax,survey,nedic2009distributed,tacle}.

A continuous-time linear consensus process over a simple connected graph $\bbb G$ can be typically modeled by a linear differential equation of the form $\dot x(t)=-Lx(t)$, where $x(t)$ is a vector in $\R^n$ and $L$ is the ``Laplacian matrix'' of $\bbb G$. 
For any simple graph $\bbb G$ with $n$ vertices, we use $D(\bbb G)$ and $A(\bbb G)$ to denote its degree matrix and adjacency matrix, respectively. Specifically, $D(\bbb G)$ is an $n\times n$ diagonal matrix whose $i$th diagonal entry equals the degree of vertex $i$, and $A(\bbb G)$ is an $n\times n$ matrix whose $ij$th entry equals 1 if $(i,j)$ is an edge in $\bbb G$ and otherwise equals 0. The Laplacian matrix of $\bbb G$ is defined as $L(\bbb G)=D(\bbb G)-A(\bbb G)$.
It is easy to see that any Laplacian matrix is symmetric and thus has a real spectrum. It is well known that $L(\bbb G)$ is positive-semidefinite, its smallest eigenvalue equals 0, and its second smallest eigenvalue, called the algebraic connectivity of $\bbb G$ and denoted as $a(\bbb G)$, is positive if and only if $\bbb G$ is connected \cite{Fiedler73}. 
It has been shown that the convergence rate of continuous-time linear consensus is determined by the algebraic connectivity, in that the larger the algebraic connectivity is, the faster the consensus can be reached \cite{reza1}.

With the preceding facts in mind, a natural and fundamental research problem is how to design network topology to achieve faster or even the fastest consensus. The problem has been studied for many years \cite{moura,ming,bayen,dai,dong,ogiwara}, to name a few. 
Notwithstanding these developments, the following question is still largely unsolved:
{\em Given a fixed number of vertices and edges, what are the optimal graphs that achieve maximal algebraic connectivity?}

\vspace{.05in}

The above question presents a challenging combinatorial optimization problem, and thus, it was only partially answered for some special cases in \cite{ogiwara}. 
Even though a powerful computer can execute such a combinatorial search, identifying the graphical properties of optimal graphs with maximal algebraic connectivity remains a mystery, not to mention the associated computational complexity.

In this paper, we propose approximating the maximal algebraic connectivity by minimizing the ``Laplacian energy'' defined as follows.

\vspace{.05in}

\begin{definition}
The Laplacian energy of a simple graph $\bbb G$ with $n$ vertices is $E(\bbb G) = \sum_{i=1}^n \lambda_i^2$, where $\lambda_i$, $i\in\{1,\ldots,n\}$ are eigenvalues of the Laplacian matrix of $\bbb G$. 
\end{definition}

\vspace{.05in}

The above concept was first proposed in \cite{Lazic06} and finds applications/connections to ordinary energy for $\pi$-electron energy in molecules \cite{Gutman78} and the first Zagreb index \cite{Gutman04}. It is worth mentioning that there have been various mathematical definitions for network energy \cite{So10}, including the earliest version of Laplacian energy \cite{Gutman06}.

\begin{figure}[!ht]
\centering
\includegraphics[width=3.43in]{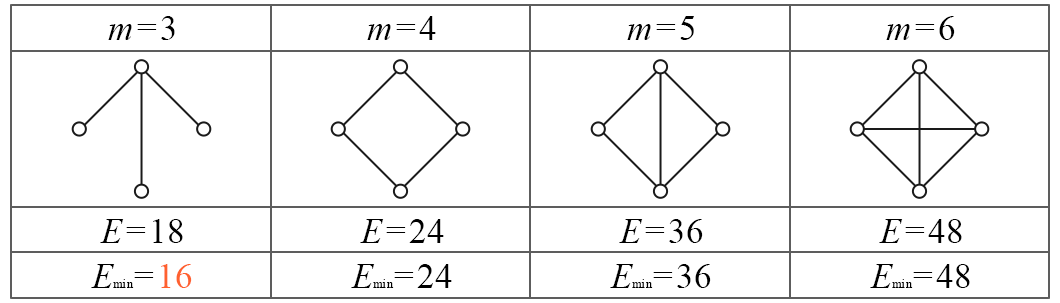} 
\caption{Maximal algebraic connectivity graphs with 4 vertices}
\label{fig:n=4}
\end{figure}

\begin{figure}[!ht]
\centering
\includegraphics[width=3.43in]{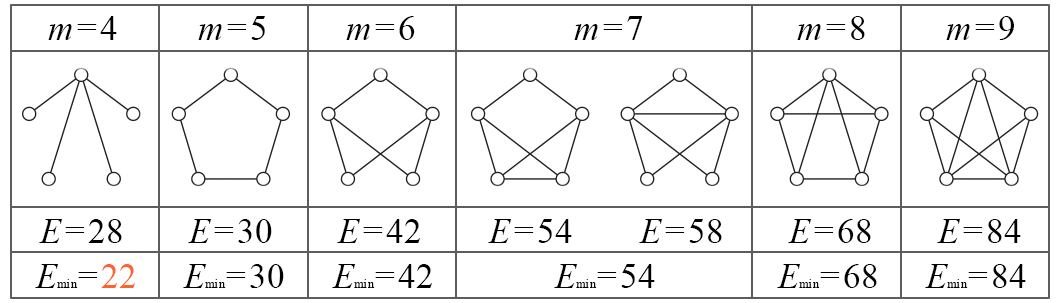} 
\caption{Maximal algebraic connectivity graphs with 5 vertices}
\label{fig:n=5}
\end{figure}

\begin{figure}[!ht]
\centering
\includegraphics[width=3.43in]{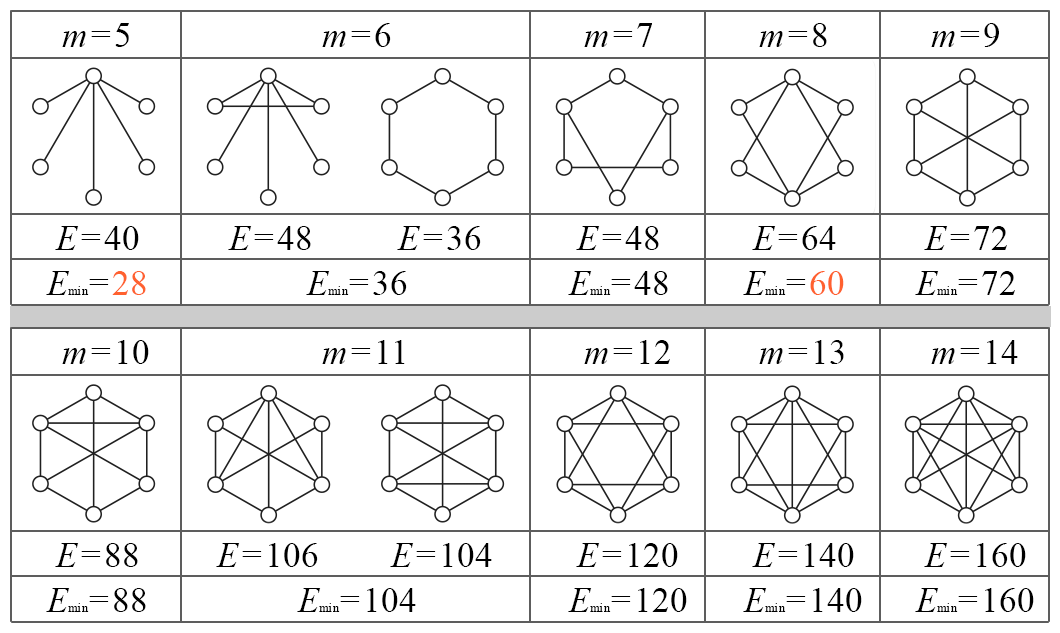} 
\caption{Maximal algebraic connectivity graphs with 6 vertices}
\label{fig:n=6}
\end{figure}

\begin{figure*}[!ht]
\centering
\includegraphics[width=7in]{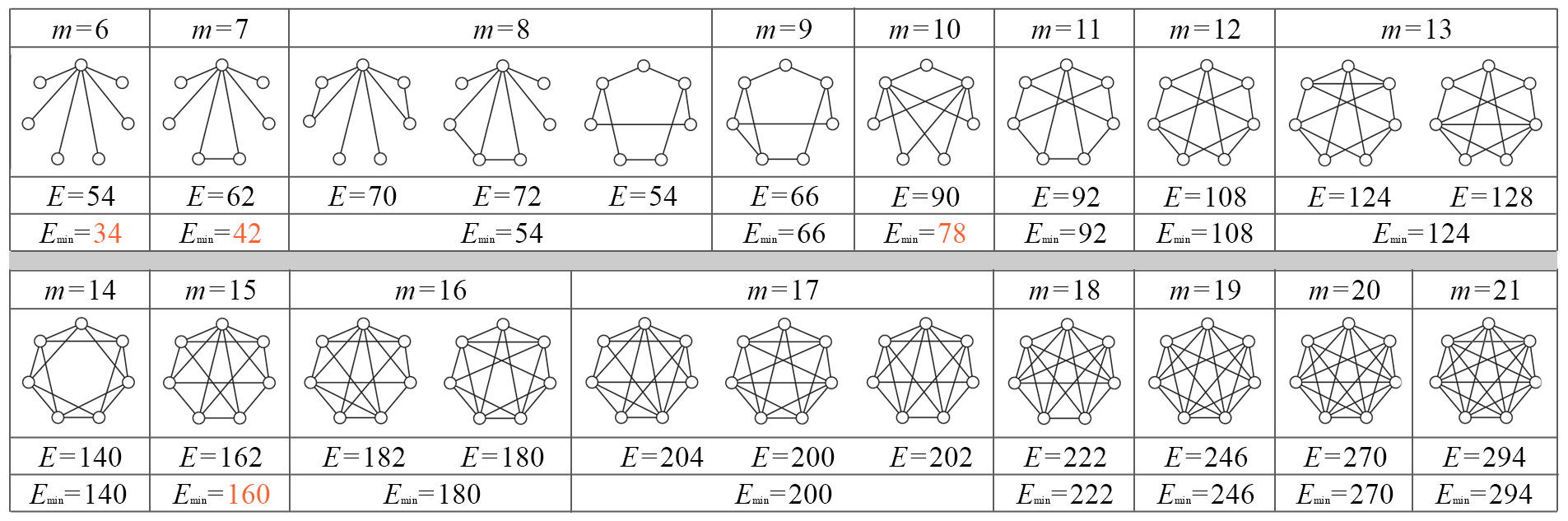} 
\caption{Maximal algebraic connectivity graphs with 7 vertices}
\label{fig:n=7}
\end{figure*}

We are motivated to appeal to the concept of Laplacian energy for designing fast/optimal consensus network topologies due to the following observations: We list all maximal algebraic connectivity graphs under the constraint of fixed numbers of vertices and edges for the cases where the vertex number $n$ ranges from 4 to 7. These are respectively given in Figures \ref{fig:n=4} through \ref{fig:n=7}.\footnote{The maximal algebraic connectivity graphs depicted in Figures \ref{fig:n=5} through \ref{fig:n=7} are sourced from \cite{ogiwara}.} We omit the case of $n=3$ as well as some complete graphs, as these graphs are unique. For $n\ge 8$, it will be very computationally expensive to go through all possible graphs. 
For each of these graphs, we list its corresponding Laplacian energy $E$, and for each pair of vertex number $n$ and edge number $m$, we list the minimal Laplacian energy $E_{{\rm min}}$ among all possible graphs. It is readily apparent that among all simple graphs with a fixed number of vertices and edges, maximal algebraic connectivity and minimal Laplacian energy coincide in most cases. The non-matching cases, highlighted in orange, are always centered in sparse cases and occasionally scattered in medium-density cases. This suggests that we may design fast consensus topologies by minimizing Laplacian energy for most scenarios. It turns out that, given a fixed number of vertices and edges, minimizing Laplacian energy is a much easier task and considerably more computationally efficient.

In this paper, we first characterize the degree distribution properties of minimal Laplacian energy graphs under the constraint of fixed numbers of vertices and edges, and then devise an algorithm to construct such connected optimal graphs (cf. Section \ref{sec:energy}). Next, we show that the minimal Laplacian energy graphs generated by the proposed algorithm exhibit strong resilience, featuring maximum vertex and edge connectivity (cf. Section \ref{sec:resilience}). Finally, we investigate the spectral properties of the Laplacian matrices of these generated minimal Laplacian energy graphs, and show that they possess large algebraic connectivity of optimal order, provided they are not sparse (cf. Section \ref{sec:fast}). Overall, we propose a computationally efficient approach to designing fast and resilient consensus topologies by minimizing Laplacian energy.

\section{Minimal Laplacian Energy}\label{sec:energy}

It has been proved in \cite{Lazic06} that $E(\bbb G) = \sum_{i=1}^n (d_i^2 + d_i)$, where $d_i$ denotes the degree of vertex $i$. It immediately implies that the Laplacian energy of a simple graph will increase after adding any additional edge. Thus, the Laplacian energy of an $n$-vertex graph achieves its maximum value, $n^2(n-1)$, when the graph is complete. Various upper and lower bounds on $E(\bbb G)$ have been established \cite{Zhou08}.
There has been an effort in the literature to identify optimal topologies that minimize the Laplacian energy for certain types of graphs. For example, among all $n$-vertex connected graphs, the Laplacian energy achieves its minimum value when the graph is the path \cite{Lazic06}. 
Another example is that among all connected graphs with chromatic number $\chi$, the Laplacian energy achieves its minimum value, $\chi^2(\chi-1)$, by the $\chi$-vertex complete graph \cite{Liu10}.
Notwithstanding these results, the following question has never been studied: 
{\em Given a fixed number of vertices and edges, what are the optimal graphs that achieve minimal Laplacian energy?}

This section solves the above open problem.
To state our main results, we use $\lfloor \cdot \rfloor$ to denote the floor function. 

\vspace{.05in}

\begin{theorem}\label{th:PNAS-symmetric}
Among all simple graphs with $n$ vertices and $m$ edges, the minimal Laplacian energy is
$(k+1)(4m-nk)$ with $k = \lfloor \frac{2m}{n} \rfloor$, which is achieved if, and only if, $n(k+1)-2m$ vertices are of degree $k$ and the remaining $2m-nk$ vertices are of degree $k+1$.
\end{theorem}

\vspace{.05in}


The theorem states that the sequence of vertex degrees, arranged in descending order, follows the following pattern:
\begin{equation}\label{eq:degrees}
(d_1, \ldots, d_n) = (\; \underbrace{k+1, \ldots, k+1}_{2m-nk}, \underbrace{k, \ldots, k}_{n(k+1)-2m})
\end{equation}
In the special case when $\frac{2m}{n}$ is an integer, all $n$ vertices are of degree $k=\frac{2m}{n}$. Thus, Theorem \ref{th:PNAS-symmetric} implies that minimal Laplacian energy graphs have an (almost) uniform degree distribution, which is intuitive from the fact that $E(\bbb G) = \sum_{i=1}^n (d_i^2 + d_i)=\sum_{i=1}^n d_i^2 + 2m$. Such a graph, whose degree difference is at most 1, is called almost regular graph. Such a graph, in which the degree difference is at most 1, is called an almost regular graph \cite{Alon84}.

It is easy to check that the total degree sum of a minimal Laplacian energy graph specified by Theorem \ref{th:PNAS-symmetric} equals $k(n(k+1)-2m)+(k+1)(2m-nk)=2m$, which is consistent with the assumption of $m$ edges. Moreover, it can be proved that such a degree distribution always admits a graph using the Erd\H{o}s-Gallai theorem \cite{Erdos60,Tripathi09}.

To prove Theorem \ref{th:PNAS-symmetric}, we need the following results. 

\vspace{.05in}

\begin{lemma}\label{lem:energy}
    (Theorem 3 in \cite{Lazic06}) For any simple graph $\bbb G$ with $n$ vertices, $E(\bbb G) = \sum_{i=1}^n (d_i^2+d_i)$.
\end{lemma}

\vspace{.05in}

\begin{lemma}\label{lem:EG}
(Erd\H{o}s-Gallai Theorem \cite{Erdos60,Tripathi09})
A nonincreasing sequence of nonnegative integers $d_1,\ldots,d_n$ constitutes the degree sequence of an $n$-vertex simple graph if, and only if, $\sum_{i=1}^n d_i$ is even and $\sum_{i=1}^j d_i \le j(j-1)+\sum_{i=j+1}^n \min\{j,d_i\}$ for all $j \in \{1,\ldots,n\}$.
\end{lemma}

\vspace{.05in}


{\bf Proof of Theorem \ref{th:PNAS-symmetric}:}
From Lemma \ref{lem:energy}, $E(\bbb G) = \sum_{i=1}^n (d_i^2+d_i)$. Since $\sum_{i=1}^n d_i = 2m$, minimizing $E(\bbb G)$ is equivalent to minimizing $\sum_{i=1}^n d_i^2$. 
Without loss of generality, assume that $d_1\ge d_2\ge \cdots \ge d_n$. 

We first minimize $\sum_{i=1}^n d_i^2$ without considering the constraint that $d_1,\ldots,d_n$ form the degree sequence of a simple graph. 
Since all $d_i$ are integers and their summation $\sum_{i=1}^n d_i$ is a fixed constant, it is easy to see that $\sum_{i=1}^n d_i^2$ is minimized by the unique sequence pattern given in \eqref{eq:degrees}, in which $k$ and $2m-nk$ are respectively the unique quotient and remainder of $2m$ divided by $n$. 
With such a degree sequence, 
$E(\bbb G) = \sum_{i=1}^n d_i^2+2m =  (k+1)^2 (2m-nk) + k^2 (n(k+1)-2m)+2m= (k+1)(4m-nk)$. 

To prove the theorem, it is now sufficient to show that, given a fixed number of vertices $n$ and edges $m$, there exists a graph with the degree sequence specified in \eqref{eq:degrees}. From Lemma \ref{lem:EG}, it is equivalent to show that (I) $\sum_{i=1}^n d_i$ is even and (II) $\sum_{i=1}^j d_i \le j(j-1)+\sum_{i=j+1}^n \min\{j,d_i\}$ for all $j \in \{1,\ldots,n\}$. 
First note that $\sum_{i=1}^n d_i=k(n(k+1)-2m)+(k+1)(2m-nk)=2m$, which implies that condition (I) holds. Next we will validate condition (II).

Since $k = \lfloor \frac{2m}{n} \rfloor$ and $m \le \frac{1}{2}n(n-1)$, it follows that $k \le n-1$. In the case when $k=n-1$, it must be true that $m=\frac{1}{2}n(n-1)$, and thus $\bbb G$ is a complete graph with all $d_i=n-1$, which satisfies \eqref{eq:degrees}. Therefore, we only need to consider $k<n-1$  from here on. To validate condition (II) for all $j\in\{1,\ldots,n\}$, we consider $j\le k$ and $j>k$ separately.  
First, in the case when $j\le k$, since $d_i\ge k$, $j\le d_i$ for all $i$. Then, $j(j-1)+\sum_{i=j+1}^n \min\{j,d_i\}=j(j-1)+(n-j)j=(n-1)j$. Meanwhile, since $d_i\le k+1$ for all $i$, $\sum_{i=1}^j d_i \le (k+1)j$. Since $k<n-1$, it follows that $\sum_{i=1}^j d_i \le (k+1)j \le (n-1)j = j(j-1)+\sum_{i=j+1}^n \min\{j,d_i\}$. Thus, condition (II) holds for all $j\le k$. 

Next, we consider the case when $j>k$. Since $d_i\le k+1$, $j\ge d_i$ for all $i$, $j(j-1)+\sum_{i=j+1}^n \min\{j,d_i\}=j(j-1)+\sum_{i=j+1}^n d_i=j(j-1)+2m-\sum_{i=1}^j d_i$, in which we used $\sum_{i=1}^n d_i=2m$ in the last equality. Then, in this case condition (II) is equivalent to 
\begin{equation}\label{eq:case2}
2\sum_{i=1}^j d_i \le j(j-1)+2m.
\end{equation}
Recalling the degree sequence given in \eqref{eq:degrees}, we will now consider the following two cases.

Case 1: $j\le 2m-nk$. From \eqref{eq:degrees}, $\sum_{i=1}^j d_i=(k+1)j$. Then, inequality \eqref{eq:case2} is equivalent to $j^2-(2k+3)j+2m \ge 0$. Let $f(x)=x^2-(2k+3)x+2m$, which is a quadratic function of $x$ and achieves its minimum at $x=k+\frac{3}{2}$. Since $j$ is an integer in the interval $[k+1,2m-nk]$, the function $f(j)=j^2-(2k+3)j+2m$ achieves its minimum at $j=k+1$, whose value is $f(k+1) = (k+1)^2 - (2k+3)(k+1) + 2m = 2m-(k^2+3k+2)$. Since $2m-nk\ge k+1$ and $k<n-1$, $2m\ge (n+1)k+1\ge (k+3)k+1$. Since $2m$ is even and $(k+3)k+1$ must be odd, $2m\ge (k+3)k+2$, which implies that $f(k+1)\ge 0$. Therefore, \eqref{eq:case2} holds in this case. 

Case 2: $j > 2m-nk$. 
From \eqref{eq:degrees}, $\sum_{i=1}^j d_i=(k+1)(2m-nk)+k(j-2m+nk)=2m-nk-kj$. Then, inequality \eqref{eq:case2} is equivalent to $j^2-(2k+1)j+2nk-2m \ge 0$. Let $f(x)=x^2-(2k+1)x+2nk-2m$, which is a quadratic function of $x$ and achieves its minimum at $x=k+\frac{1}{2}$. 
Recall that $j>k$ is under consideration. Then, $j$ is an integer in the interval $[\max\{k+1,2m-nk+1\},n]$. We will consider two scenarios, $k> 2m-nk$ and $k\le 2m-nk$, separately. 

First assume $k> 2m-nk$, which is equivalent to $k>\frac{2m}{n+1}$. Then, $j$ is an integer in the interval $[k+1,n]$, and thus the function $f(j)=j^2-(2k+1)j+2nk-2m$ over this interval achieves its minimum at $f(k+1)=-k^2+(2n-1)k-2m$. Treat this expression as a quadratic function of $k$ which achieves its maximum at $k=n-\frac{1}{2}$. Since $k$ is an integer in the interval $(\frac{2m}{n+1},n-1)$, the function over this interval satisfies $f(k+1)> f(\frac{2m}{n+1}+1)=\frac{2m}{(n+1)^2}[n(n-1)-2-2m]$. Since $k<n-1$, implying $n(n-1)>2m$, and both $n(n-1)$ and $2m$ are even, $n(n-1)-2\ge 2m$. Therefore, $f(k+1)> 0$, and thus \eqref{eq:case2} holds in this scenario. 

Next assume $k \le 2m-nk$, which is equivalent to $k \le \frac{2m}{n+1}$. Then, $j$ is an integer in the interval $[2m-nk+1,n]$, and thus the function $f(j)=j^2-(2k+1)j+2nk-2m$ over this interval achieves its minimum at $f(2m-nk+1)=n(n+2)k^2-(4mn+4m-n+2)k+4m^2$. 
Note that if $k=0$, this minimum value equals $4m^2\ge 0$, and thus \eqref{eq:case2} holds. Hence, in the remainder of the proof, we need only consider $k\ge 1$. Also note that with $k\ge 1$, if $n=2$, for which $m$ may equal 0 or 1, this minimum value is always positive, and thus \eqref{eq:case2} holds. Hence, we also need only consider $n\ge 3$ in the remainder of the proof. 
Treat this expression as a quadratic function of $k$ which, without any restriction on $k$, achieves its minimum value of  $-16m^2+(8n^2-8n-16)m - (n-2)^2$ when $k=\frac{4mn+4m-n+2}{2n(n+2)}$. But it has been assumed that $k \le \frac{2m}{n+1}$. This minimum value is achievable only if $\frac{4mn+4m-n+2}{2n(n+2)} \le \frac{2m}{n+1}$, which is equivalent to $m\le \frac{1}{4}(n^2-n-2)$. 
Also treat this minimum value as a quadratic function of $m$ which achieves its maximum at $m=\frac{1}{4}(n^2-n-2)$. 
Since $k=\lfloor\frac{2m}{n}\rfloor\ge 1$, it follows that $m\ge \frac{n}{2}$. Then, the quadratic function of $m$, $g(m)=-16m^2+(8n^2-8n-16)m - (n-2)^2$ over the interval $[\frac{n}{2}, \frac{1}{4}(n^2-n-2)]$ is always no smaller than $g(\frac{n}{2})=4n^3 - 9n^2 - 4n - 4$, which can be straightforwardly verified to be greater than 0 for all $n\ge 3$. This implies that \eqref{eq:case2} holds. On the other hand, if $\frac{4mn+4m-n+2}{2n(n+2)} > \frac{2m}{n+1}$, which is equivalent to $m> \frac{1}{4}(n^2-n-2)$. Then, the quadratic function of $k$, $f(2m-nk+1)=n(n+2)k^2-(4mn+4m-n+2)k+4m^2$ over the interval $[1,\frac{2m}{n+1}]$, achieves its minimum when $k=\frac{2m}{n+1}$, with the corresponding minimum value being $\frac{2m}{(n+1)^2)}[n(n-1)-2-2m]$, which was proven to be positive at the end of the last paragraph. Therefore, \eqref{eq:case2} always holds in this scenario.
\hfill $\qed$


\vspace{.05in}

\begin{figure}[!ht]
\centering
\includegraphics{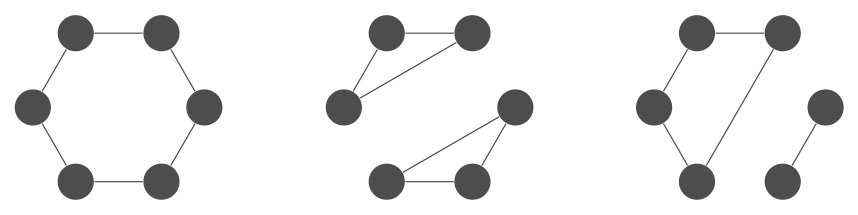} 
\caption{All minimal Laplacian energy graphs with 6 vertices and 6 edges}
\label{fig:optimal}
\end{figure}

Theorem \ref{th:PNAS-symmetric} provides a simple graphical condition, dependent solely on the degree distribution, that characterizes minimal Laplacian energy graphs. The following example shows that such a degree distribution condition does not guarantee connectivity.
Consider all simple graphs with 6 vertices and 6 edges (i.e., $n=m=6$), Theorem \ref{th:PNAS-symmetric} identifies three minimal Laplacian energy graphs, as illustrated in Figure \ref{fig:optimal}.

In many network applications, such as distributed control and optimization \cite{survey,arc_distopt}, connected graphs are often desired. The following theorem states that a connected optimal graph with the same minimal Laplacian energy always exists. 

\vspace{.05in}

\begin{theorem}\label{th:PNAS-symmetric-connected}
Among all connected simple graphs with $n$ vertices and $m$ edges, the minimal Laplacian energy is
$(k+1)(4m-nk)$ with $k = \lfloor \frac{2m}{n} \rfloor$, which is achieved if, and only if, $n(k+1)-2m$ vertices are of degree $k$ and the remaining $2m-nk$ vertices are of degree $k+1$.
\end{theorem}

\vspace{.05in}

{\bf Proof of Theorem \ref{th:PNAS-symmetric-connected}:}
In light of Theorem \ref{th:PNAS-symmetric}, it is sufficient to show that there exists a connected graph whose degree sequence satisfies \eqref{eq:degrees}. Such a graph must indeed exist, as it follows directly from Theorem \ref{thm:algorithmbasic}.
\hfill $\qed$

\vspace{.05in}

The above theorem implicitly assumes that $m\ge n-1$, otherwise there is no such connected graph.
In a special case when $m=n-1$, all connected graphs are trees. It is not hard to see that the following result is a direct consequence of Theorem \ref{th:PNAS-symmetric-connected}.

\vspace{.05in}

\begin{corollary}\label{coro:tree}
Among all simple trees with $n$ vertices, the minimal Laplacian energy is $6n-8$, which is achieved by the path.
\end{corollary}

\vspace{.05in}

\begin{figure}[!ht]
\centering
\includegraphics{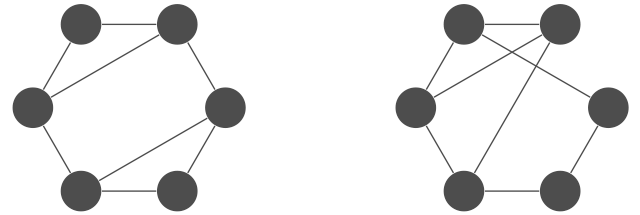} 
\caption{Two minimal Laplacian energy connected graphs with 6 vertices and 8 edges}
\label{fig:multiple-optimal-connected}
\end{figure}

For general pairs of $n$ and $m$ with $m > n-1$, the optimal connected graph may not be unique. For example, for the case when $n=6$ and $m=8$, two connected minimal Laplacian energy graphs are illustrated in Figure \ref{fig:multiple-optimal-connected}.

We next present the following algorithm, which provides a procedure to construct a minimal Laplacian energy graph that is connected and satisfies the degree distribution specified by Theorem~\ref{th:PNAS-symmetric-connected}.

\vspace{.05in}

\textit{Algorithm 1:} Given $n$ and $m$ with $m\ge n-1>0$, without loss of generality, label $n$ vertices from 1 to $n$. Set  $k = \lfloor \frac{2m}{n} \rfloor$, which implies that $nk\le 2m<n(k+1)$. 

\vspace{.05in}

\noindent\hspace{.05in}
{\bf Case 1:} The integer $k$ is even.
\begin{enumerate}
    \item[(1)] If $2m=nk$, for each $i\in\{1,\ldots,n\}$, connect vertex $i$ with each of those vertices whose indices are $(i+j)\bmod n$, $j\in\{\pm 1,\ldots,\pm\frac{k}{2}\}$.
    
    \item[(2)] If $2m=nk+l$ where $l\in[1,n)$ is an even integer, first construct the graph as done in Case 1 (1), and then for each $i\in\{1,\ldots,\frac{l}{2}\}$, connect vertex~$i$ and vertex $i+\lfloor \frac{n}{2} \rfloor$.
\end{enumerate}

\noindent\hspace{.05in}
{\bf Case 2:} The integer $k$ is odd.

\begin{enumerate}
    \item[(1)] If $n$ is even and $2m=nk$, for each $i\in\{1,\ldots,n\}$, connect vertex $i$ with each of those vertices whose indices are $(i+j)\bmod n$, $j\in\{\pm 1,\ldots,\pm\frac{k-1}{2},\frac{n}{2}\}$.
    
    \item[(2)] If $n$ is even and $2m=nk+l$ where $l\in[1,n)$ is an even integer, first construct the graph as done in Case 2 (1), and then for each $i\in\{1,\ldots,\frac{l}{2}\}$, connect vertex~$i$ and vertex $i+\frac{n-2}{2}$.

    \item[(3)] If $n$ is odd and $2m=nk+1$, first for each $i\in\{1,\ldots,n\}$, connect vertex $i$ with each of those vertices whose indices are $(i+j)\bmod n$, $j\in\{\pm 1,\ldots,\pm\frac{k-1}{2}\}$, and then for each $i\in\{1,\ldots,\frac{n+1}{2}\}$, connect vertex $i$ and vertex $i+\frac{n-1}{2}$.

    \item[(4)] If $n$ is odd and $2m=nk+1+l$ where $l\in[1,n-1)$ is an even integer, first construct the graph as done in Case 2 (3), and then for each $i\in\{\frac{n+3}{2},\ldots,\frac{n+1+l}{2}\}$, connect vertex~$i$ and vertex $(i+\frac{n-1}{2}) \bmod n$.
\end{enumerate}

\vspace{.05in}

\begin{theorem}\label{thm:algorithmbasic}
    Algorithm 1 constructs a connected simple graph with $n(k+1)-2m$ vertices of degree $k$ and $2m-nk$ vertices of degree $k+1$. 
\end{theorem}

\vspace{.05in}

{\bf Proof of Theorem \ref{thm:algorithmbasic}:}
First of all, any graph generated by Algorithm 1 must be connected, as each case within the algorithm contains a cycle with a vertex sequence $(1,2,\ldots,n,1)$.
Next, it is straightforward to verify that the graphs generated by each case of the algorithm satisfy the degree sequence specified by \eqref{eq:degrees}.
\hfill $\qed$

\vspace{.05in}

\begin{figure}[!ht]
\centering
\includegraphics[width=2.3in]{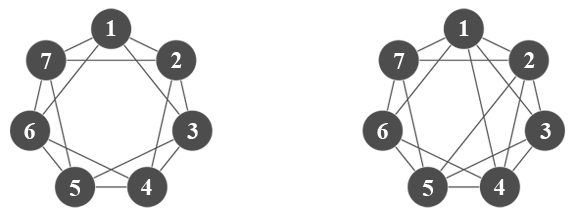} 
\caption{Left is the graph with 7 vertices and 14 edges, generated by Algorithm 1 following Case 1 (1); right is the graph with 7 vertices and 16 edges, generated by Algorithm~1 following Case 1 (2).
}
\label{fig:case1}
\end{figure}

\begin{figure}[!ht]
\centering
\includegraphics[width=2.3in]{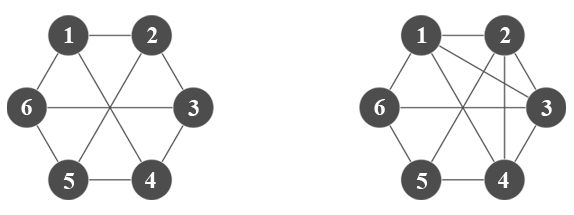} 
\caption{Left is the graph with 6 vertices and 9 edges, generated by Algorithm 1 following Case 2 (1); right is the graph with 6 vertices and 11 edges, generated by Algorithm~1 following Case 2 (2).
}
\label{fig:case2-(1)(2)}
\end{figure}

\begin{figure}[!ht]
\centering
\includegraphics[width=2.3in]{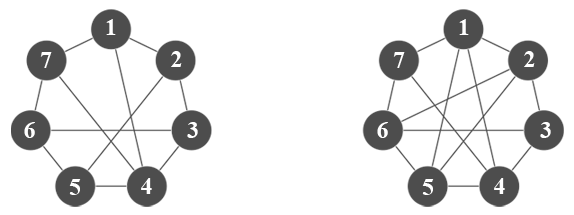} 
\caption{Left is the graph with 7 vertices and 11 edges, generated by Algorithm 1 following Case 2 (3); right is the graph with 7 vertices and 13 edges, generated by Algorithm~1 following Case 2 (4).
}
\label{fig:case2-(3)(4)}
\end{figure}

It can be straightforwardly checked that in the case when $n=m=6$, Algorithm 1 will follow Case 1 (1) and construct the 6-vertex path, which is consistent with Corollary \ref{coro:tree}. We further present six tailored examples, each corresponding to a distinct case outlined in Algorithm 1, as depicted in Figures \ref{fig:case1} to \ref{fig:case2-(3)(4)}. These examples collectively validate Theorem \ref{thm:algorithmbasic}.

\section{Connectivity Resilience}\label{sec:resilience}

The graphs generated by Algorithm~1 exhibit ``optimal'' connectivity properties. To see this, we introduce two well-known connectivity concepts in graph theory: vertex connectivity, denoted as $v(\bbb G)$, which is defined as the minimum number of vertices whose removal would disconnect graph $\bbb G$, and edge connectivity, denoted as $e(\bbb G)$, which is defined as the minimum number of edges whose removal would disconnect graph $\bbb G$.
For the complete graph with $n$ vertices, it is obvious that its edge connectivity equals $n-1$. However, there is no subset of vertices whose removal disconnects the complete graph. It is conventional to set its vertex connectivity as $n-1$ \cite[page 149]{West00}.

\vspace{.05in}

\begin{theorem}\label{th:algorithm}
    Let $\bbb G$ be the graph generated by Algorithm 1 with $n$ vertices and $m$ edges. Then, $v(\bbb G) = e(\bbb G)  =\lfloor \frac{2m}{n} \rfloor$. 
\end{theorem}

\vspace{.05in}


\vspace{.05in}

The theorem implies that the graphs constructed by Algorithm 1 always have the maximum vertex and edge connectivity. To see this, consider any graph $\bbb G$ with $n$ vertices and $m$ edges. Its minimum degree $\delta (\bbb G)$ is at most $\lfloor \frac{2m}{n} \rfloor$. Since $e(\bbb G)\le \delta (\bbb G)$ by definition and $v(\bbb G)\le e(\bbb G)$ \cite[Theorem 5]{whitney}, it follows that $v(\bbb G)\le e(\bbb G)\le \delta (\bbb G)\le \lfloor \frac{2m}{n} \rfloor$.  

\begin{figure}[!ht]
\centering
\includegraphics[width=2.3in]{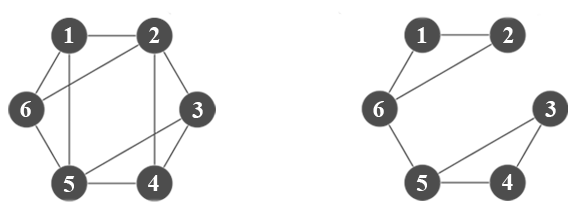} 
\caption{Two almost regular graphs
}
\label{fig:almost}
\end{figure}

Recall that all graphs generated by Algorithm 1 are almost regular graphs with the degree sequence specified in \eqref{eq:degrees}. In general, almost regular graphs do not necessarily have $v(\bbb G) = e(\bbb G) = \lfloor \frac{2m}{n} \rfloor$. 
To see this, consider two examples in Figure \ref{fig:almost}. The left almost regular graph has $n=6$ vertices and $m=10$ edges, but its vertex connectivity is 2 (by removing vertex 2 and vertex 5), which is smaller than $\lfloor \frac{2m}{n} \rfloor=3$. The right almost regular graph has $n=6$ vertices and $m=7$ edges, but its vertex connectivity is 1 (by removing vertex 6), so is its edge connectivity (by removing the edge between vertex 5 and vertex 6), both being smaller than $\lfloor \frac{2m}{n} \rfloor=2$.


\vspace{.05in}


To prove Theorem \ref{th:algorithm}, we need the following notation and lemma. Let $\scr U$ be a vertex subset of graph $\bbb G$. We use $\bbb G\setminus\scr U$ to denote the graph resulting from removing all vertices in $\scr U$, along with the edges connecting any vertex in $\scr U$, from~$\bbb G$. 


\vspace{.05in}

\begin{lemma}\label{lm:alg}
Let $\bbb G$ be a simple graph with $n$ vertices and $p<n$ be a positive integer. Suppose that for any pair of distinct vertices $i,j\in\{\1,2,\ldots,n\}$ with $|i-j|\le (p \bmod n)$, $(i,j)$ is an edge in $\bbb G$. Let $\scr U$ be any vertex subset of $\bbb G$ and $u,v\in\{\1,2,\ldots,n\}$ be any two vertices in $\bbb G$ such that $u<v$ and $u,v\notin\scr U$. If there are no $p$ vertices in $\scr U$ whose indices are consecutive integers in the interval $(u,v)$, then there exists a path between $u$ and $v$ in $\bbb G\setminus\scr U$.
\end{lemma}


\vspace{.05in}

{\bf Proof of Lemma \ref{lm:alg}:}
We explicitly construct a such path between $u$ and $v$. Consider a sequence of vertex indices: $u_0=u$ and for any $i \ge 1$, 
\begin{itemize}
    \item If $u_i < v-p$, define $u_{i+1}$ to be the largest index in $[u_i+1, u_i+p]$ such that vertex $u_{i+1}$ is not in $\scr U$. Such $u_{i+1}$ always exists because there are no $p$ vertices in $\scr U$ whose indices are consecutive integers in $(u,v)$.
    \item Otherwise, the sequence ends.
\end{itemize}
By construction, vertices $u_i$ and $u_{i+1}$ are adjacent. Let $u_j$ be the final term of this sequence. Then $v-p \le u_j \le v$.
\begin{itemize}
    \item If $u_j = v$, we have found the path $u=u_0 \to u_1 \to \dots \to u_j=v$.
    \item Otherwise, vertex $u_j$ is adjacent to $v$, so we have found the path $u=u_0 \to u_1 \to \dots \to u_j \to v$.
\end{itemize}
Either way, there is a path between $u$ and $v$.
\hfill $\qed$

\vspace{.05in}

{\bf Proof of Theorem \ref{th:algorithm}:}
From the preceding discussion,
$v(\bbb G)\le e(\bbb G)\le k=\lfloor \frac{2m}{n} \rfloor$. Thus, to prove the theorem, it is sufficient to show $v(\bbb G) \ge k$. 

Graphs constructed in Case 1 (2) are formed by adding edges to graphs in Case 1 (1); graphs constructed in Case 2 (2) are formed by adding edges to graphs in Case 2 (1); graphs constructed in Case 2 (4) are formed by adding edges to graphs in Case 2 (3). Since vertex connectivity is non-decreasing when we add edges, it suffices to check $v(\bbb G) \ge k$ for graphs from Case 1 (1), Case 2 (1), and Case 2 (3) only.

To prove $v(\bbb G) \ge k$, we will equivalently show that for any subset $\scr U$ of $k-1$ vertices, $\bbb G \setminus \scr U$ is connected. Say that a vertex is deleted if it is in $\scr U$. Let $u,v$ be two arbitrary non-deleted vertices, where $u < v$. It suffices to show that there is a path between $u$ and $v$ in $\bbb G \setminus \scr U$.

\vspace{0.05in}

\textbf{Case 1 (1):} From the algorithm description, there is an edge between any two vertices whose indices differ by at most $k/2$. Since $k-1$ vertices are deleted, we have either
\begin{itemize}
    \item There are no $k/2$ deleted vertices whose indices are consecutive integers in $(u,v)$.
    \item Or there are no $k/2$ deleted vertices whose indices are consecutive integers in $\{1, \dots, u-1\} \cup \{v+1, \dots, n\}$.
\end{itemize}
In either scenario, Lemma \ref{lm:alg} implies that there is a path between $u$ and $v$ in $\bbb G \setminus \scr U$.

\vspace{0.05in}

\textbf{Case 2 (1):} From the algorithm description, there is an edge between any two vertices whose indices differ by at most $(k-1)/2$ and an edge between any vertex $i$ and vertex $i+n/2$. If either
\begin{itemize}
    \item There are no $(k-1)/2$ deleted vertices whose indices are consecutive integers in $(u,v)$.
    \item Or there are no $(k-1)/2$ deleted vertices whose indices are consecutive integers in $\{1, \dots, u-1\} \cup \{v+1, \dots, n\}$.
\end{itemize}
then Lemma \ref{lm:alg} implies that there is a path between $u$ and $v$ in $\bbb G \setminus \scr U$.

Thus, it remains to consider the case in which there are $(k-1)/2$ deleted vertices whose indices are consecutive integers in $(u,v)$ and $(k-1)/2$ deleted vertices whose indices are consecutive integers in $\{1, \dots, u-1\} \cup \{v+1, \dots, n\}$. Let the former set of vertex indices be $\scr A$ and the latter be $\scr B$. Since $|\scr U| = k-1$, we see that $\scr U = \scr A \cup \scr B$. Let $\scr B = \{b, \dots, b+(k-3)/2\}$, where indices are taken mod $n$.

\vspace{.05in}

If $v = u+n/2$, we directly obtain the path $u \to v$, which is what we wanted to show. Thus, suppose $v \neq u+n/2$. We will only consider the case $v < u+n/2$ because the case $v > u+n/2$ is symmetric.

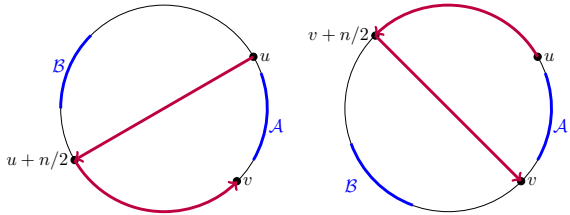
\begin{figure}[h!]
    \centering
\resizebox{3in}{!}{\begin{tikzpicture}
\draw[color=black](0,0) circle (2);
\filldraw[black] (1.73,1) circle (2pt) node[anchor=west]{$u$};
\filldraw[black] (-1.73,-1) circle (2pt) node[anchor=east]{$u+n/2$};
\filldraw[black] (1.41,-1.41) circle (2pt) node[anchor=west]{$v$};
\draw[blue, ultra thick] (1.73,-1) arc (-30:20:2);
\node (a) at (2.17,-0.34) {\color{blue} $\scr A$};
\draw[blue, ultra thick] (-2,0) arc (180:135:2);
\node (b) at (-2.04,0.82) {\color{blue} $\scr B$};
\draw[purple, ultra thick, ->] (1.73,1) -- (-1.73,-1);
\draw[purple, ultra thick, ->] (-1.73,-1) arc (210:315:2);

\draw[color=black](5.5,0) circle (2);
\filldraw[black] (5.5+1.73,1) circle (2pt) node[anchor=west]{$u$};
\filldraw[black] (5.5+1.41,-1.41) circle (2pt) node[anchor=west]{$v$};
\filldraw[black] (5.5-1.41,1.41) circle (2pt) node[anchor=east]{$v+n/2$};
\draw[blue, ultra thick] (5.5+1.73,-1) arc (-30:20:2);
\node (a) at (5.5+2.17,-0.34) {\color{blue} $\scr A$};
\draw[blue, ultra thick] (3.62,-0.68) arc (200:250:2);
\node (b) at (3.62,-1.5) {\color{blue} $\scr B$};
\draw[purple, ultra thick, ->] (5.5-1.41,1.41) -- (5.5+1.41,-1.41);
\draw[purple, ultra thick, ->] (7.23,1) arc (30:135:2);
\end{tikzpicture}}
\caption{Illustration of Case 2 (1)}
\label{fig:thm4}
\end{figure}

Since $\scr A \subseteq (u,v)$, we have $v-u > (k-1)/2$, so $(v+n/2)-(u+n/2) > (k-1)/2$. Then the length of the interval $(u+n/2, v+n/2)$ is greater than $|\scr B|$.
\begin{itemize}
    \item If $b > u+n/2$, then $\scr B$ does not contain $u+n/2$, so consider the path $u \to u+n/2 \to u+n/2-1 \to \dots \to v$. See the left illustration in Figure \ref{fig:thm4} for this case.
    \item Otherwise, $b \le u+n/2$. Since the length of the interval $(u+n/2, v+n/2)$ is greater than $|\scr B|$, we see that $\scr B$ does not contain $v+n/2$. Consider the path $u \to u-1 \to \dots \to v+n/2+1 \to v+n/2 \to v$. See the right illustration in Figure \ref{fig:thm4} for this case.
\end{itemize}
In both scenarios, we found a path between $u$ and $v$ in $\bbb G \setminus \scr U$.

\vspace{0.05in}

\textbf{Case 2 (3):} From the algorithm description, there is an edge between any two vertices whose indices differ by at most $(k-1)/2$ and an edge between $i$ and $i+(n-1)/2$ for any $1 \le i \le \frac{n+1}{2}$.

Using the same argument as in Case 2 (1), Lemma \ref{lm:alg} solves all scenarios, except the one where there are $(k-1)/2$ deleted vertices whose indices are consecutive integers in $(u,v)$ and $(k-1)/2$ deleted vertices whose indices are consecutive integers in $\{1, \dots, u-1\} \cup \{v+1, \dots, n\}$. Let the former set of vertex indices be $\scr A$ and the latter be $\scr B$. Since $|\scr U| = k-1$, we see that $\scr U = \scr A \cup \scr B$. Let $\scr B = \{b, \dots, b+(k-3)/2\}$, where indices are taken mod $n$.

\vspace{.05in}

The condition that there is an edge between $i$ and $i+(n-1)/2$ for any $1 \le i \le \frac{n+1}{2}$ translates to:
\begin{itemize}
    \item Vertices with index at most $(n-1)/2$ are adjacent to $i+(n-1)/2$
    \item Vertices with index at least $(n+3)/2$ are adjacent to $i+(n+1)/2$
    \item Vertex $(n+1)/2$ is adjacent to both $(n+1)/2+(n-1)/2$ and $(n+1)/2+(n+1)/2$.
\end{itemize}

Let $u$ be adjacent to $u+(n + \varepsilon_u)/2$, and $v$ be adjacent to $v+(n + \varepsilon_v)/2$, where $\varepsilon_u, \varepsilon_v \in \{-1,1\}$. We claim that the length of the interval $(u+(n + \varepsilon_u)/2, v+(n + \varepsilon_u)/2)$ is greater than $(k-1)/2$.

Since $\scr A \subseteq (u,v)$, we have $v-u \ge (k+1)/2$, so the length of the interval is at least $(k-1)/2$. For the length to equal $(k-1)/2$, vertex $u$ must be adjacent to $u+(n+1)/2$, and vertex $v$ must be adjacent to $v+(n-1)/2$. If $u,v \neq (n+1)/2$, we immediately reach a contradiction because $u < v$. If $u=(n+1)/2$, change $\varepsilon_u$ to $-1$. If $v$ equals $(n+1)/2$, change $\varepsilon_v$ to $1$. Thus, we can always ensure that the length of the interval is greater than $(k-1)/2$.

\begin{itemize}
    \item If $b > u+(n + \varepsilon_u)/2$, then $\scr B$ does not contain $u+(n + \varepsilon_u)/2$, so consider the path $u \to u+(n + \varepsilon_u)/2 \to u+(n + \varepsilon_u)/2-1 \to \dots \to v$.
    \item Otherwise, $b < u+(n + \varepsilon_u)/2$. Since the length of the interval $(u+(n + \varepsilon_u)/2, v+(n + \varepsilon_u)/2)$ is greater than $|\scr B|$, we see that $\scr B$ does not contain $v+(n + \varepsilon_v)/2$, so consider the path $u \to u-1 \to \dots \to v+(n + \varepsilon_v)/2+1 \to v+(n + \varepsilon_v)/2 \to v$.
\end{itemize}
In both scenarios, we found a path between $u$ and $v$ in $\bbb G \setminus \scr U$.

\vspace{.05in}

In all cases, we have found a path between any two vertices in $\bbb G \setminus \scr U$, so $\bbb G \setminus \scr U$ is connected. Thus, $v(\bbb G) = k$, which completes the proof.
\hfill $\qed$


\section{Fast Consensus}\label{sec:fast}

In this section, we study the algebraic connectivity of the minimal Laplacian energy graphs generated by Algorithm 1. We will show that the generated ``dense'' graphs possess large algebraic connectivity, while the generated ``sparse'' graphs do not. This finding is consistent with the observations from the figures in the introduction.

Among all non-complete graphs with $n$ vertices and $m$ edges, it is known that $a(\bbb G) \le v(\bbb G)$ \cite[Theorem 4.1]{Fiedler73}. From the preceding discussion, it follows that $a(\bbb G) \le \lfloor \frac{2m}{n} \rfloor$ for all non-complete graphs. Theorem \ref{th:algorithm-a} gives a lower bound on algebraic connectivity of graphs constructed by Algorithm 1.



\begin{theorem}\label{th:algorithm-a}
Let $\bbb G$ be the graph generated by Algorithm~1 with $n$ vertices and $m\ge n$ edges. Then, 
\eq{a(\bbb G) \ge \bar k - \frac{\sin (\bar k\pi / n)}{\sin (\pi/n)},\label{eq:thm5}}
where $k=\lfloor \frac{2m}{n} \rfloor$ and $\bar k = 2\lfloor \frac{k}{2} \rfloor+1$, with equality holding if the graph is constructed in Case 1 (1) and Case 2 (1). 
\end{theorem}

\vspace{.05in}

Note that $\bar k \le k+1 \le \frac{2m}{n} +1\le  n$.
We will use this fact without special mention in the sequel.

To prove Theorem \ref{th:algorithm-a}, we need the following concept and results. A circulant matrix is a square matrix in which all rows are composed of the same entries and each row is rotated one entry to the right relative to the preceding row. The spectrum of any circulant matrix can be completely determined by its first row entries, as specified in the following lemma: 

\vspace{.05in}

\begin{lemma}\label{lm:part2-1}
(Theorem 6 in \cite{Kra12})
If $C$ is an $n\times n$ circulant matrix whose first row entries are $c_0, c_1,\ldots, c_{n-1}$, then its $n$ eigenvalues are $\lambda_i = \sum_{p=0}^{n-1} c_p e^{\frac{j2pi \pi}{n}}$, $i\in\{0,1,\ldots,n-1\}$, where $j$ is  the imaginary unit. 
\end{lemma}






\vspace{.05in}

\begin{lemma}\label{lm:part2-2}
For any integers $n\ge 2$ and $2 \le k \le n-2$,
$$\max_{i\in\{1,2,\ldots,n-1\}} 2\sum_{p=1}^{\lfloor \frac{k}{2} \rfloor} \cos \bigg(\frac{2pi\pi }{n} \bigg) = \frac{\sin (\bar k\pi / n)}{\sin (\pi/n)}-1,$$
where $\bar k = 2\lfloor \frac{k}{2} \rfloor+1$, and the maximum is achieved if, and only if, $i=1$ or $i=n-1$.
\end{lemma}


\vspace{.05in}

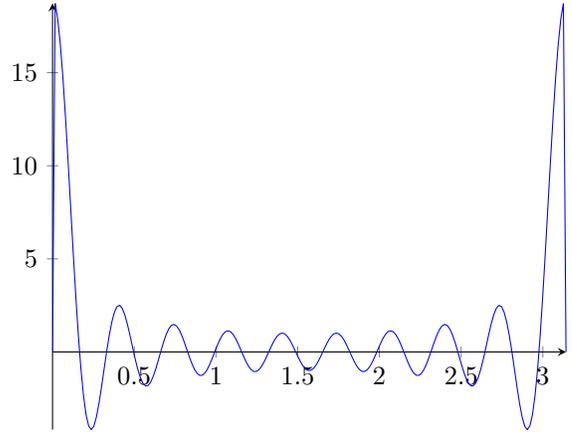
\begin{figure}[!ht]
    \centering
    \resizebox{3in}{!}{
        \begin{tikzpicture}
        \begin{axis}[axis lines = center]
          \addplot[domain=0:pi,samples=200,blue]{sin(19*deg(x)) / sin(deg(x))};
        \end{axis}
        \end{tikzpicture}}
    \caption{Plot of $f(x)=\frac{\sin(\bar k x)}{\sin(x)}$ for $\bar{k}=19$}
    \label{fig:f(x)}
\end{figure}

{\bf Proof of Lemma \ref{lm:part2-2}:}
From the angle addition and subtraction formulae, $2\cos a\sin b = \sin(a+b)-\sin(a-b)$ for any $a,b\in\R$. Let $a=2pi\pi/n$ and $b=i\pi/n$. Since $\sin (i\pi/n)>0$ for all integers $i\in\{1,2,\ldots,n-1\}$, it follows that for all integers $i\in\{1,2,\ldots,n-1\}$,
\[ 2\cos \bigg(\frac{2pi\pi}{n} \bigg) = \frac{\sin ((2p+1) i\pi /n ) - \sin ((2p-1) i\pi /n )}{\sin (i\pi /n )}. \]
Summing this relation over index $p$ from $1$ to $\lfloor \frac{k}{2} \rfloor$, 
\[ 2\sum_{p=1}^{\lfloor \frac{k}{2} \rfloor} \cos \bigg(\frac{2pi\pi}{n} \bigg) = \frac{\sin (\bar k i\pi /n)}{\sin (i\pi /n)}-1. \]
It remains to identify the optimal integers $i$ over the interval $[1,n-1]$ that maximize $\frac{\sin (\bar k i\pi /n)}{\sin (i\pi /n)}$. 
To this end, define a function $f(x) = \frac{\sin(\bar k x)}{\sin(x)}$. Figure \ref{fig:f(x)} provides an example plot.
First, we check that $f(x)$ is symmetric around $\pi/2$. Since $\bar{k} = 2\lfloor k/2 \rfloor+1$ is odd, 
\[ f(\pi-x)= \frac{\sin(\bar{k} (\pi-x))}{\sin(\pi-x)} = \frac{\sin(\pi-\bar{k}x)}{\sin(\pi-x)} = \frac{\sin(\bar{k}x)}{\sin(x)} = f(x). \]
Thus, it suffices to show $f(\pi/n) \ge f(i\pi/n)$ for all integers $i \in (0,n/2]$.
Second, we show that in the interval $(0,\pi)$, $f(x)$ attains its maximum of $\bar k$ as $x$ approaches $0$ or $\pi$. Clearly, $\lim_{x\to 0} f(x) = \lim_{x\to \pi} f(x) = \bar k$. Meanwhile,
\begin{align*}
|f(x)| &= \bigg|\frac{e^{j\bar k x} - e^{-j\bar k x}}{e^{jx} - e^{-jx}} \bigg| \\
&= \big|e^{j(\bar k-1)x}+e^{j(\bar k-3)x}+\dots+e^{j(-\bar k+1)x}\big| \\
&\le \big|e^{j(\bar k-1)x}\big|+\big|e^{j(\bar k-3)x}\big|+\dots+\big|e^{j(-\bar k+1)x}\big| \le \bar k.
\end{align*}
Finally, observe that $f(x)$ has roots at $\pi/\bar{k}, 2\pi/\bar{k}, \dots, (\bar{k}-1)\pi/\bar{k}$ and has a maximum/minimum in each interval $(p\pi/\bar{k}, (p+1)\pi/\bar{k})$, where $0 \le p \le \bar{k}-1$. Since $k \le n-2$ by assumption, we have $\bar{k} \le n-1$, so $\pi/n$ is always in the leftmost interval and is closest to $x$-coordinate of the global maximum at 0. Hence, $f(\pi/n) \ge f(i\pi/n)$ for all integers $i \in (0,n/2]$.
\hfill $\qed$


\vspace{.05in}

{\bf Proof of Theorem \ref{th:algorithm-a}:}
Note that for any simple graph with $n\ge 2$ vertices and $m\ge n$ edges, $k=\lfloor \frac{2m}{n}\rfloor \ge 2$. Also, if $k = n-1$, the graph is complete, so $a(\bbb G) = n$ \cite{Fiedler73} and inequality \eqref{eq:thm5} holds. Thus, we only consider $k \le n-2$ in the rest of the proof. Since the algebraic connectivity of a non-complete graph will not decrease after adding an edge \cite[Corollary 3.2]{Fiedler73}, it suffices to prove \eqref{eq:thm5} for Case 1 (1), Case 2 (1), and Case 2 (3).

\vspace{.05in}

\textbf{Case 1 (1):} From the algorithm description, it is straightforward to verify that the Laplacian matrix of the generated graph is a circulant matrix, with its first row entries being
$$k,\; \underbrace{-1, \ldots, -1}_{k/2}, \;\underbrace{0, \ldots, 0}_{n-k-1}, \;\underbrace{-1, \ldots, -1}_{k/2}$$
From Lemma \ref{lm:part2-1}, its $n$ eigenvalues are
\begin{align*}
\lambda_i &= k - \sum_{p=1}^{k/2} e^{\frac{j 2pi \pi}{n}} - \sum_{p=n-k/2}^{n-1} e^{\frac{j 2pi \pi}{n}}  \\
&= k - \sum_{p=1}^{k/2} \left(e^{\frac{j 2pi \pi}{n}} + e^{\frac{j 2(n-p)i \pi}{n}} \right) \\
&= k - 2\sum_{p=1}^{k/2} \cos \bigg(\frac{2pi \pi}{n} \bigg), \;\;\;i\in\{0,1,\ldots,n-1\}.
\end{align*}
It is easy to see that $\lambda_0 = 0$, which is the smallest eigenvalue. Since the generated graph is connected, all other eigenvalues are positive. 
Since in Case 1, $k$ is even, and thus $\bar k= k+1$. From Lemma \ref{lm:part2-2}, the maximum among $2\sum_{p=1}^{k/2} \cos \left(\frac{2pi \pi}{n} \right)$, $i\in\{1,2,\ldots,n-1\}$ is $\frac{\sin ((k+1)\pi / n)}{\sin (\pi/n)} - 1$. Thus, the second smallest eigenvalue 
$a(\bbb G) = k+1 - \frac{\sin ((k+1)\pi / n)}{\sin (\pi/n)}$. 

\vspace{.05in}

\textbf{Case 2 (1):}
From the algorithm description, it is straightforward to verify that the Laplacian matrix of the generated graph is a circulant matrix, with its first row entries being
$$k,\; \underbrace{-1, \ldots, -1}_{(k-1)/2}, \;\underbrace{0, \ldots, 0}_{(n-k-1)/2}, \; -1, \;\underbrace{0, \ldots, 0}_{(n-k-1)/2}, \;\underbrace{-1, \ldots, -1}_{(k-1)/2}$$
From Lemma \ref{lm:part2-1}, its $n$ eigenvalues are
\begin{align*}
\lambda_i &= k - \sum_{p=1}^{(k-1)/2} e^{\frac{j 2pi \pi}{n}} -1 - \sum_{p=n-(k-1)/2}^{n-1} e^{\frac{j 2pi \pi}{n}}  \\
&= k - 1 - \sum_{p=1}^{(k-1)/2} \left(e^{\frac{j 2pi \pi}{n}} + e^{\frac{j 2(n-p)i \pi}{n}} \right) \\
&= k - 1 - \sum_{p=1}^{(k-1)/2} \left(e^{\frac{j 2pi \pi}{n}} + e^{-\frac{j 2pi \pi}{n}} \right) \\
&= k-1 - 2\sum_{p=1}^{(k-1)/2} \cos \bigg(\frac{2pi \pi}{n} \bigg),\;\;\;i\in\{0,1,\ldots,n-1\}.
\end{align*}
Since in Case 2, $k$ is odd, and thus $\bar k= k$. Using the same argument as in Case 1 (1), the second smallest eigenvalue 
$ a(\bbb G) = k - \frac{\sin (k\pi / n)}{\sin (\pi/n)}$.

\vspace{.05in}

\textbf{Case 2 (3):}  
In this case, the Laplacian matrix of the generated graph $\bbb G$ is not a circulant matrix. Consider the spanning subgraph of $\bbb G$, denoted as $\bbb H$, with an edge set defined such that for each pair of $i\in\{1,\ldots,n\}$ and $j \in \{\pm 1, \dots, \pm \frac{k-1}{2} \}$, there is an edge between vertex $i$ and vertex $(i+j) \bmod n$.
It is straightforward to verify that $\bbb H$ is connected and its Laplacian matrix is a circulant matrix, with its first row entries being
$$k,\; \underbrace{-1, \ldots, -1}_{(k-1)/2}, \;\underbrace{0, \ldots, 0}_{n-k}, \;\underbrace{-1, \ldots, -1}_{(k-1)/2}$$
Using the same argument as in Case 1 (1), the second smallest eigenvalue 
$ a(\bbb H) = k+1 - \frac{\sin (k\pi / n)}{\sin (\pi/n)}$. 
Since $\bbb H$ is a spanning subgraph of $\bbb G$, $a(\bbb G)\ge  k+1 - \frac{\sin (k\pi / n)}{\sin (\pi/n)}$.
\hfill$\qed$

\vspace{.1in}

From Theorem \ref{th:algorithm-a}, the graphs constructed by Algorithm 1 in Case 1 (1) and Case 2 (1) have an explicit algebraic connectivity expression, $a(\bbb G) = \bar k - \frac{\sin (\bar k\pi / n)}{\sin (\pi/n)}$, which can be bounded as follows:

\vspace{.05in}

\begin{lemma}\label{th:case1-1}
For any graph $\bbb G$ generated by Algorithm 1 in Case 1 (1) or Case 2 (1),
$\frac{\pi^2 (0.5\bar{k}^3-\bar{k})}{6n^2 - \pi^2}< a(\bbb G)< \frac{\bar{k}^3\pi^2}{6n^2}$, where $\bar k = 2\lfloor \frac{k}{2} \rfloor+1$ and $k=\lfloor \frac{2m}{n}\rfloor$.
\end{lemma}

\vspace{.05in}

{\bf Proof of Lemma \ref{th:case1-1}:}
Using the Taylor series and basic calculus, it is straightforward to show that $x-x^3/6 < \sin x < x-x^3/6+x^5/120$ for all $x>0$. 
Applying the upper bound in this inequality to $\sin(\bar{k}\pi/n)$ leads to $\sin(\bar{k}\pi/n) < \frac{\bar{k}\pi}{n} - \frac{\bar{k}^3\pi^3}{6n^3} + \frac{\bar{k}^5\pi^5}{120n^5}$, and applying the lower bound to $\sin(\pi/n)$ leads to $\sin(\pi/n) > \frac{\pi}{n} - \frac{\pi^3}{6n^3}$. Then, 
\begin{align*}
a(\bbb G) &> \bar{k} - \frac{\frac{\bar{k}\pi}{n} - \frac{\bar{k}^3\pi^3}{6n^3} + \frac{\bar{k}^5\pi^5}{120n^5}}{\frac{\pi}{n} - \frac{\pi^3}{6n^3}} = \frac{(\bar{k}^3-\bar{k}) n^2\pi^2 - \frac{1}{20}\bar{k}^5\pi^4 }{6n^4 - n^2\pi^2 } \\
&= \frac{(\bar{k}^3-\bar{k})\pi^2}{6n^2 - \pi^2} - \frac{\frac{1}{20} \bar{k}^5\pi^4}{6n^4 - n^2\pi^2}.
\end{align*}
Since $\bar{k} \le n$, it follows that $\frac{\bar{k}^5}{6n^4 - \pi^2n^2} \le \frac{\bar{k}^3}{6n^2 - \pi^2}$. Thus,
\begin{align*}
a(\bbb G) &> \frac{(\bar{k}^3-\bar{k})\pi^2}{6n^2 - \pi^2} - \frac{\frac{1}{20} \bar{k}^3\pi^4}{6n^2 - \pi^2} 
> \frac{\pi^2 (0.5\bar{k}^3-\bar{k})}{6n^2 - \pi^2}.
\end{align*}
Next we apply $x-x^3/6 < \sin x$ to $\sin(\bar{k}\pi/n)$, which leads to $\sin(\bar{k}\pi/n) > \frac{\bar{k}\pi}{n}- \frac{\bar{k}^3\pi^3}{6n^3}$. With this and the fact $\sin(\pi/n) < \pi/n$, it follows that $a(\bbb G)< \frac{\bar{k}^3\pi^2}{6n^2}$. 
\hfill $\qed$

\vspace{.1in}

More can be said. Let $\lfloor \frac{2m}{n} \rfloor \le \sqrt[3]{6n^2/\pi^2} -1$. Then, $\bar{k} \le k+1 \le \sqrt[3]{6n^2/\pi^2}$. From Lemma \ref{th:case1-1}, $a(\bbb G) < 1$. We have thus proved the following:

\vspace{.05in}

\begin{corollary}\label{coro:sparse}
Let $\bbb G$ be any graph generated by Algorithm 1 in Case 1 (1) or Case 2 (1). If $k=\lfloor \frac{2m}{n} \rfloor \le \sqrt[3]{6n^2/\pi^2} -1$, then $a(\bbb G) < 1$.
\end{corollary}

\vspace{.05in}

Let us agree to call a graph with $n$ vertices and $m$ edges sparse if its average degree $\frac{2m}{n}$ is much smaller than $O(n)$, and dense if $\frac{2m}{n} = O(n)$.

Among all connected graphs with $n$ vertices and $m$ edges, it is known that $a(\bbb G) \ge 2e(\bbb G) (1-\cos(\pi/n))$ \cite[Theorem 4.3]{Fiedler73}. Since $e(\bbb G) \le n-1$, it is easy to see that this lower bound of $a(\bbb G)$ is strictly less than 1 if $n \ge 9$. Corollary \ref{coro:sparse} implies that when the minimal Laplacian energy graph generated by Algorithm 1 in Case 1 (1) or Case 2 (1) is sparse, its algebraic connectivity is small. This suggests that small/minimal Laplacian energy and large/maximal algebraic connectivity do not match for sparse graphs. This observation is not surprising; for instance, in the special case of tree graphs where $m=n-1$, the maximal algebraic connectivity graph is the star, while the minimal Laplacian energy graph is the path, which are opposites. The maximal algebraic connectivity graphs for some special sparse cases were theoretically identified in \cite[Theorems 1, 2, 4]{ogiwara}.

In contrast to sparse graphs, the following result shows that generated dense graphs have large algebraic connectivity.

\vspace{.05in}

\begin{corollary}\label{coro:newdense}
Let $\bbb G$ be the graph constructed by Algorithm 1 with $n$ vertices and $m$ edges. If $k=\lfloor \frac{2m}{n} \rfloor\ge n+1-\sqrt{2n-3}$, then $a(\bbb G)\ge k - 2\sqrt{k-1}$.
\end{corollary}

\vspace{.05in}

Among all connected non-complete graphs with $n$ vertices and $m$ edges, it is known that $a(\bbb G) \le v(\bbb G)$ \cite[Theorem 4.1]{Fiedler73}. 
From the discussion in Section \ref{sec:resilience}, $a(\bbb G) \le v(\bbb G) \le k$. Note that $n+1-\sqrt{2n-3}=O(n)$. Thus, Corollary \ref{coro:newdense} implies that when the minimal Laplacian energy graph generated by Algorithm 1 is dense, its algebraic connectivity is large. 
In fact, it is known that $a(\bbb G) \le k - 2\sqrt{k-1} + O(\log_k n)^{-1}$ for all connected almost regular graphs with $n$ vertices and $m$ edges \cite[Theorem 1]{Alon91}. Therefore, the dense graphs generated by Algorithm 1 exhibit nearly optimal algebraic connectivity.

\vspace{.05in}

\begin{figure}[!ht]
    \centering
    \resizebox{3in}{!}{
        \begin{tikzpicture}
        \begin{axis}[axis lines = center]
          \addplot[domain=0:25,samples=200,red]{sin(deg(x*pi/24)) / sin(deg(pi/24))-(2*x-1)^0.5};
          \addplot[domain=0:25,samples=200,blue]{-1.143*x+20.589};
          \filldraw[black] (23,-5.708) circle (1pt) node[anchor=east]{\footnotesize $(x_0,f(x_0))$};
          \filldraw[black] (18.007,0) circle (1pt) node[anchor=south]{\footnotesize $(x_1,0)$};
        \end{axis}
        \end{tikzpicture}}
    \caption{Plot of Newton's method on $f(x)$ for $n=24$. The red curve is $f(x)$, and the blue line is the tangent.}
    \label{fig:corollary2}
\end{figure}
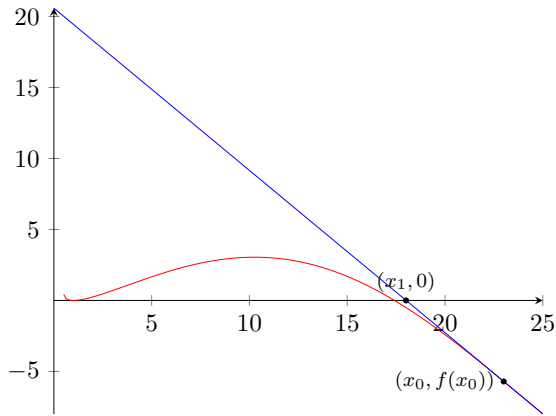

{\bf Proof of Corollary \ref{coro:newdense}:} 
From Theorem \ref{th:algorithm-a}, to prove the corollary, it is sufficient to show
$\bar k - \frac{\sin (\bar k\pi / n)}{\sin (\pi/n)}\ge k - \sqrt{2k-1}$. 
In the case when $k$ is odd, $\bar k = k$, and the target inequality simplifies to $\frac{\sin(k\pi/n)}{\sin(\pi/n)} \le \sqrt{2k-1}$.
In the case when $k$ is even, $\bar k = k+1$, and the target inequality simplifies to $\frac{\sin((k+1)\pi/n)}{\sin(\pi/n)} - 1 \le \sqrt{2k-1}$. Since $k$ is an integer in the interval $[n+1-\sqrt{2n-3}, n-1]$ and it is easy to verify that $n/2 \le n+1-\sqrt{2n-3}$, it follows that $n/2 \le n+1-\sqrt{2n-3} \le k \le n-1$.
Then, $\pi/2 \le k\pi/n < (k+1)\pi/n \le \pi$. Since the sine function is decreasing on the interval $[\pi/2, \pi]$,  $\sin(k\pi/n) > \sin((k+1)\pi/n)$. Thus, it is sufficient to prove $\frac{\sin(k\pi/n)}{\sin(\pi/n)} \le \sqrt{2k-1}$ for all even $k$ in the interval $[n+1 - \sqrt{2n-3}, n-1]$. 
To this end, we define the following function: 
\[ f(x) = \frac{\sin (x\pi / n)}{\sin (\pi/n)} - \sqrt{2x-1}, \;\;\; x \in [n+1 - \sqrt{2n-3}, n-1]\]
We then apply Newton's method to obtain a bound on the largest root of $f(x)$ in the interval $[n+1 - \sqrt{2n-3}, n-1]$. Figure \ref{fig:corollary2} is a visualization for $n=24$. Set the initial estimate $x_0 = n-1$. The first estimation $x_1$ is the intersection of the tangent line at $(x_0,f(x_0))$ with the $x$-axis:\begin{align*}
x_1 &= x_0 - \frac{f(x_0)}{f'(x_0)} = n-1-\frac{\sqrt{2n-3}-1}{\frac{\pi}{n} \cot(\frac{\pi}{n}) + (2n-3)^{-\frac{1}{2}}} \\
&\le n-1-\frac{\sqrt{2n-3}-1}{1 + (2n-3)^{-\frac{1}{2}}} \le n+1-\sqrt{2n-3}.
\end{align*}
Since $f(x)$ is concave downward on the interval $[n+1 - \sqrt{2n-3}, n-1]$, the actual largest root must be strictly less than $x_1$. Hence, for all $x \ge n+1-\sqrt{2n-3}$, $f(x)$ must be nonpositive, which completes the proof.
\hfill $\qed$

The graphs constructed by Algorithm 1 in Case 1 (1), which actually belong to the so-called regular lattices. A simple graph with $n\ge 3$ vertices is called a $d$-regular lattice, with $d$ being an even integer in the interval $[2,n-1]$, if each vertex $i$ is adjacent to each of those vertices whose indices are $(i+j)\bmod n$, $j\in\{\pm 1,\ldots,\pm\frac{d}{2}\}$ \cite{Watts98}.  It is easy to see that any graph with $n\ge 3$ vertices generated by Algorithm 1 in Case 1 (1) is a $k$-regular lattice.
Lemma \ref{th:case1-1} immediately implies that the algebraic connectivity of a $k$-regular lattice is of the order $O(k^3/n^2)$.


Regular lattices are closely related to Watts-Strogatz small-world networks, which are generated by randomly rewiring edges in a regular lattice. The rewiring procedure involves iterating through each edge, and with probability $p$, one endpoint is moved to a new vertex chosen randomly from the lattice. Double edges and self-loops are not allowed in this process, so small-world networks are simple graphs \cite{Watts98}.


The work of \cite{Olfati05} defines the algebraic connectivity gain, $\lambda_2(p)/\lambda_2(0)$,  as the algebraic connectivity of the small-world network formed by rewiring with probability $p$ divided by the algebraic connectivity of the regular lattice \cite[Definition 1]{Olfati05}. By running simulations with $k \approx \log(n)$, the paper conjectures that the maximum $\lambda_2(p)/\lambda_2(0)$ is on the order of $O(n)$ \cite[Observation (ii), page 4]{Olfati05}, which implies that small-world networks can reach consensus significantly faster than regular lattices.
Lemma \ref{th:case1-1} implies that $\lambda_2(0)$ is on the order of $O(k^3/n^2)$. Thus, if one can show that $\lambda_2(p)$ is on the order of $O(k^3/n)$, then the conjecture in \cite{Olfati05} would be mathematically verified. Whether small-world networks can achieve algebraic connectivity of order $O(k^3/n)$ has so far eluded us, but Lemma \ref{th:case1-1} provides a helpful first step in addressing this question.







\section{Conclusion}


This paper proposes a novel approach to designing fast consensus topologies by minimizing Laplacian energy, marking the first step in this direction.
Although both sparse and dense graphs have been analyzed, and the findings are consistent with the observations in the introduction, the scattered non-matching cases for medium-dense graphs (see Figures \ref{fig:n=6} and \ref{fig:n=7}) have not been addressed. These graphs belong to complete bipartite graphs. A simple graph is called bipartite if its vertices can be partitioned into two classes so that every edge has endpoints in different classes. The complete bipartite graph $K_{p,n-p}$ is the bipartite graph with $p$ vertices in one class, $n-p$ vertices in the other class, and all $p(n-p)$ edges between vertices of different classes \cite[page 17]{Diestel}. It is worth mentioning that maximal algebraic connectivity graphs were identified as complete bipartite graphs for some special medium-dense graphs \cite[Theorem 3, Table I]{ogiwara}. Understanding these observations is a direction for future research.




\bibliographystyle{unsrt}
\bibliography{susie,jicareer,consensus,push}

\end{document}